\documentclass[12pt,leqno]{amsart}
\usepackage{amscd,times,fullpage}
\usepackage[utf8]{inputenc}
\usepackage{amsmath}
\usepackage{amsthm}
\usepackage{mathtools}
\usepackage{amsfonts}
\usepackage{comment}
\usepackage{graphicx}
\usepackage{enumerate}
\usepackage{url}
\usepackage{tikz}
\usetikzlibrary{cd}
\usepackage[hidelinks]{hyperref}
\usepackage{cleveref}
\usepackage{makecell}
\usepackage[cmtip,all]{xy}
\newcommand{\longsquiggly}{\xymatrix{{}\ar@{~>}[r]&{}}}

\hypersetup {
	pdfpagemode = {UseNone},
	pdftitle = {K\"ahler manifolds are not that formal},
	pdfauthor = {Giovanni Placini, Jonas Stelzig, Leopold Zoller},
	pdflang = {en-UK}
} 

\newtheorem{thmintro}{Theorem}

\newtheorem{prop}{Proposition}
\newtheorem{lem}[prop]{Lemma}

\newtheorem{definition}[prop]{Definition}
\newtheorem{thm}[prop]{Theorem}
\newtheorem{cor}[prop]{Corollary}
\theoremstyle{remark}
\newtheorem{remark}[prop]{Remark}
\newtheorem{rem}[prop]{Remark}

\newtheoremstyle{colon}%
{}
{}
{\itshape}
{}
{\bfseries}
{:}
{ }
{}

\theoremstyle{colon}

\newcommand{\C}{\mathbb{C}}

\newcommand{\CP}{\mathbb{CP}}

\newcommand{\Z}{\mathbb{Z}}

\newcommand{\im}{\operatorname{im}}

\newcommand{\del}{\partial}

\newcommand{\delbar}{{\bar\partial}}

\newcommand{\Oh}{\mathcal{O}}

\usepackage[all]{xy}

\newcommand{\cD}{\mathcal{D}}

\makeatother

\usepackage{multirow}

\title{Nontrivial Massey products on compact K\"ahler manifolds}

\author{Giovanni Placini}
\address{Dipartimento di Matematica e Informatica, Universit\`a degli Studi di Cagliari, via Ospedale 72, 09124 Cagliari, Italy}
\email{giovanni.placini@unica.it}

\author{Jonas Stelzig}
\address{Mathematisches Institut der LMU M\"unchen, Theresienstra{\ss}e 39, 80333 M\"unchen, Germany}
\email{jonas.stelzig@math.lmu.de}

\author{Leopold Zoller}
\address{Facultat de Matemàtiques i Informàtica, Universitat de Barcelona, Gran Via de les Corts Catalanes, 585, 08007 Barcelona}
\email{leopold.zoller@ub.edu}

\begin{document}
\thanks{The first named author is supported by INdAM and GNSAGA and by ISI-HOMOS. He acknowledges financial support by PNRR e.INS (CUP F53C22000430001, codice MUR ECS00000038). The visits of the authors to LMU and UniCA were funded by the programs MGR, GOACT, and Visiting Professors of UniCA}
\begin{abstract}
By a result of Deligne, Griffiths, Morgan, and Sullivan, compact K\"ahler manifolds are rationally formal. In contrast to this, we show that the pluripotential quasi-isomorphism type of the bigraded, bidifferential algebra of forms on a compact K\"ahler manifold generally contains more information than the de Rham cohomology algebra with its real Hodge structure. More precisely, on any compact K\"ahler manifold fibering over a Riemann surface of genus at least two, there is a nontrivial $ABC$-Massey product. Furthermore, starting from dimension three, there are simply connected projective manifolds with a nonzero $ABC$-Massey product of three divisor classes. In particular, compact K\"ahler manifolds are generally not formal in the sense of pluripotential homotopy theory. 
\end{abstract}
\maketitle

\section{Introduction}

The quasi-isomorphism class of the commutative differential graded algebra (cdga) of differential forms provides a much finer algebraic invariant of manifolds than the de Rham cohomology algebra. For instance, one may compute from it topological higher operations such as Massey products and, under assumptions of finite type and nilpotent $\pi_1$, the ranks of the higher homotopy groups \cite{Sullivan}.

One of the early highlights of Sullivan's approach to rational homotopy theory is the theorem of \cite{DGMS} that any compact K\"ahler manifold is (de Rham) formal: The cdga of forms can be connected by a chain of quasi-isomorphisms to the de Rham cohomology algebra, considered as a cdga with trivial differential. In other words, the quasi-isomorphism class of the differential forms is encoded in the de Rham cohomology and hence the former does not contain more information than the latter. In particular, on compact Kähler manifolds all Massey products need to vanish.

For complex manifolds, the splitting $d=\del+\delbar$ equips the complex-valued differential forms with the enriched structure of a commutative bigraded, bidifferential algebra (cbba). 
We consider cbbas up to pluripotential quasi-isomorphisms, cf. \cite{StePHT}, meaning maps of cbbas which induce an isomorphism in the Bott-Chern and Aeppli cohomologies
\[H_{BC}=\frac{\ker\del\cap\ker\delbar}{\im\del\delbar}\quad \text{ and } \quad H_A=\frac{\ker{\del\delbar}}{\im\del+\im\delbar}.\] 
This notion of quasi-isomorphism is in a certain sense universal \cite[Thm C]{StePHT}. In particular, it is stronger than requiring the map to induce an isomorphism in de Rham or Dolbeault cohomology.  A complex manifold is called strongly formal if its cbba of forms is connected by a chain of pluripotential quasi-isomorphisms to a cbba with trivial differentials \cite{MilSt_bigrform}. Hermitian symmetric spaces, all compact K\"ahler solvmanifolds, and all compact K\"ahler manifolds of dimension $\geq 2$ with cohomology of complete intersection type (the K\"ahler class generates all cohomology below middle degree) are strongly formal  \cite{MilSt_bigrform}, \cite{SfTo_BCK}, \cite{StePHT}.

The goal of this paper is to show that, nevertheless, compact K\"ahler manifolds are generally not strongly formal. I.e., the pluripotential quasi-isomorphism class of the cbba of forms contains more information than the de Rham cohomology algebra with its real Hodge structure. In particular, bigraded variants of the established tools of rational homotopy theory provide a new way of studying compact K\"ahler manifolds. 

To prove this, we consider a secondary variant of Massey products, called Aeppli-Bott-Chern-Massey products ($ABC$-Massey products), introduced in \cite{AnToBCform}.\footnote{In many cases, in particular those appearing in \Cref{Thm: CP3}, they can be seen to compute higher products introduced in \cite{Den_higher}. This will be explained in forthcoming work.} They are defined like ordinary Massey products, replacing $d$ by $i\del\delbar$. We recall their precise definition in \Cref{sec: prelim}. Schematically, they yield a partially defined map of bidegree $(-1,-1)$
\[
\xymatrix{
\langle\cdot,\cdot,\cdot\rangle_{ABC}: H_{BC}^{\otimes 3}(X)\ar@{-->}[r]& H_A(X)}.
\]

For compact K\"ahler manifolds, the natural maps $H_{BC}(X)\to H_{dR}(X)\to H_A(X)$ are isomorphisms. Therefore, one might view $\langle\cdot,\cdot,\cdot\rangle_{ABC}$ as a partially defined operation of degree $-2$ in de Rham cohomology, compatible with the real Hodge structures. As with de Rham formality and ordinary Massey products, on any strongly formal manifold, all $ABC$-Massey products are trivial \cite{MilSt_bigrform}.
\begin{thmintro}\label{Thm: RS} On any closed Riemann surface of genus at least two
there is a nontrivial $ABC$-Massey product.
\end{thmintro}
As a consequence, on any compact K\"ahler manifold with a surjective map to a curve of genus at least two there is a nontrivial $ABC$-Massey product, \Cref{cor: fibrations}. On the other hand, simply connectedness and nilpotency play a large role in  rational homotopy theory. It is thus natural to ask whether simply connected compact K\"ahler manifolds which are not strongly formal exist as well. Simply connected compact complex surfaces have cohomology of complete intersection type and hence are strongly formal. Therefore, the minimal dimension where this could occur is $3$. In that regard our main result is:
\begin{thmintro}\label{Thm: CP3}
    There is a finite sequence of blow-ups in points and lines of $\CP^3$ such that on the resulting space there is a nontrivial $ABC$-Massey product.
\end{thmintro}
As a consequence of Theorem \ref{Thm: CP3}, every complex manifold of dimension $\geq 4$ is bimeromorphic to one on which there is a nontrivial $ABC$-Massey product, \Cref{Cor: blowups nonform}.
The $ABC$-Massey product in \Cref{Thm: RS} takes combinations of holomorphic and antiholomorphic $1$-forms as input. In \Cref{Thm: CP3}, one takes divisor classes as inputs and the output is related to the cross ratio of four points on a line in $\CP^3$. In particular, in the latter example varying the positions of the centers of the blow-ups yields diffeomorphic manifolds which are, as complex manifolds, distinguished by their nontrivial $ABC$-Massey product, \Cref{biholotypes}. This illustrates that the variation of the pluripotential homotopy type with respect to the underlying formal rational homotopy type provides interesting new functions on moduli spaces.

\subsection*{Related work}
The question of formality of K\"ahler manifolds in pluripotential homotopy theory lies in between several formality and nonformality results: 

Adapting the proof of \cite{DGMS}, the authors of \cite{NT} show that in Dolbeault homotopy theory, where one considers forms up to Dolbeault quasi-isomorphism, compact K\"ahler manifolds are still formal. The Dolbeault homotopy type can be obtained from the pluripotential homotopy type. On general compact complex manifolds there are many higher operations encoded in the pluripotential homotopy type which are, generally speaking, not encoded by Dolbeault or de Rham homotopy theory \cite{Den_higher}, \cite{MilSt_bigrform}, \cite{TaDiss}. None of these were known to be nontrivial on compact K\"ahler manifolds. An important indication that such operations may also be interesting in the K\"ahler case was given by \cite{SfTo_DBC} who exhibited a (non K\"ahler) manifold satisfying the $\del\delbar$-Lemma with a nontrivial $ABC$-Massey product. This is relevant as the proof of formality of compact K\"ahler manifolds in the de Rham and Dolbeault settings depends only on the $\del\delbar$-Lemma. Finally, taking into account information beyond forms, \cite{CaCleMo_pi3} exhibits an example which shows that the rational Mixed Hodge Structure on homotopy groups is  generally not a formal consequence of the rational Hodge Structure on the cohomology ring. Their argument uses intermediate Jacobians and hence depends on the rational structure. It appears to us that the $ABC$-Massey products considered in this paper provide different information than Abel-Jacobi type invariants, cf. \Cref{rem: no odd}, even when considering infinitesimal ones, see \Cref{rem: rigid}.

\subsection*{Acknowledgements}  Part of the work was carried out during the visit of the first named author to LMU Munich and a visit of the second and third named authors to the University of Cagliari. We thank both institutions for their hospitality. We further thank J.~Cirici, C.~Deninger, D.~Kotschick, A.~Milivojevic, J.~Morgan, M.~Paulsen, D.~Phong, G.~Quick, S.~Rollenske, C.~Voisin for conversations about the subject. We also thank an anonymous referee for their comments, which in particular improved the proof of \Cref{Thm: RS}.

\section{Preliminaries and notation}\label{sec: prelim}
\subsection{Gradings} Given any bigraded vector space $H^{\bullet,\bullet}$, a single superscript denotes the associated singly graded space with total grading $H^\bullet=\bigoplus_{p+q=\bullet}H^{p,q}$. When we talk about the total space, we omit the grading completely $H=\bigoplus_{p,q}H^{p,q}$ and for an element $h\in H$ we denote its (bi)degree by $|h|$.

\subsection{Forms, currents and fundamental classes} For any manifold $X$, we write $A^{\bullet}(X)$ for its graded-commutative differential graded algebra of $\C$-valued smooth differential forms and $\mathcal{D}(X)$ for the space of currents, i.e. $\cD^k(X)$ is the topological dual of the space of compactly supported $\left( \dim(X)-k\right)$-forms. We write $d$ for the exterior differential and we will usually omit the wedge product $\alpha\wedge\beta=\alpha\beta$. 

Let us now assume $X$ is a complex manifold. Then $A(X)$ (and $\cD(X)$) carry a bigrading by type, refining the single grading and making it into a graded-commutative, bigraded, bidifferential algebra (resp. a bigraded, bidifferential module over $A_X$). We write $\del$ and $\delbar$ for the components of bidegree $(1,0)$, respectively $(0,1)$, of the exterior derivative $d=\del+\delbar$. For $X$ compact, the map $A(X)\to \cD(X), \omega\mapsto \int (-1)^{|\omega|-1}\omega\wedge\_$ is a pluripotential quasi-isomorphism, see \cite{Serre}, \cite[Lem.~1.27]{StePHT}. In particular one can (and we will) represent cohomology classes both by forms and by currents. As the most relevant example for us, any complex submanifold $Z\subseteq X$ of complex codimension $p$, or any formal linear combination of such, gives rise to a current of integration $\delta_Z:=(\omega\mapsto \int_Z\omega|_Z)\in \cD^{p,p}(X)$. This current is $\del$ and $\delbar$-closed and therefore represents a class in $H_{BC}^{p,p}(X)$, which we denote by $[Z]:=[\delta_Z]$. It is called the fundamental class, or Thom class, of $Z$. While this current representative is very natural, we will need also smooth form representatives for $[Z]$. These always exist, and if $Z$ is of complex codimension $1$, i.e. a divisor, the Lelong-Poincar\'e formula produces particularly concrete ones as follows: Let $\Oh(Z)$ be the associated line bundle, $s$ a meromorphic section of $\Oh(Z)$ with divisor of zero and poles $(s)=Z$ and $h$ a Hermitian metric on $\Oh(Z)$. Further, let $c_1^h(Z)$ be the form representative for $c_1(\Oh(Z))\in H_{BC}^{1,1}(X)$ built as the curvature of the Chern connection with respect to $h$. Concretely, if locally $h=e^{-\varphi}\|~\|^2$, then $c_1^h(Z)=\frac{i}{2\pi}\del\delbar\varphi$. Now $\log|s|_h$ is a locally integrable function and hence integration against it defines a current in $\cD^{0,0}(X)$. Then:
    \begin{prop}[Lelong-Poincar\'e formula]\label{prop:LLP} With the above setup one has
        $$\frac i \pi \del\delbar \log |s|_h=\delta_Z-c_1^h(Z).$$
    \end{prop}
We refer to \cite{Demailly_book} for a proof and further background on currents.

\begin{rem} \label{rem: product currents}\,
    Partial product structures exist on currents, which extend those on forms but also include transverse intersection of fundamental currents: For instance, if $A=\sum a_i A_i,B=\sum b_i B_i$ are analytic cycles all of whose components intersect transversally, then one may define $\delta_A\cdot \delta_B:=\delta_{A\cap B}$ where $A\cap B=\sum a_ib_j A_i\cap B_j$ is the linear combination of the set-theoretic intersections and if $\eta$ is a locally integrable (possibly not smooth) form which is still locally integrable when restricted to $A$, one sets $\eta\cdot\delta_A:=(\omega\mapsto \int_A\eta\omega)$. These products are still graded-commutative and the differentials satisfy a Leibniz rule, whenever all sides of the relevant equations are defined.  
\end{rem}

\subsection{Massey products}

 We now introduce the secondary version of classical Massey products this paper is concerned with. The definition of classical Massey products is recalled at the end of this section.

\begin{definition}[\cite{AnToBCform}]\label{def: ABCMP}
    Let $X$ be a complex manifold and $[\alpha]\in H_{BC}^{k_1}(X)$, $[\beta]\in H_{BC}^{k_2}(X)$, $[\gamma]\in H_{BC}^{k_3}(X)$ be s.t. $[\alpha][\beta]=[\beta][\gamma]=0$. The triple Aeppli-Bott-Chern Massey product of these classes is the coset
    \[
        \langle\alpha,\beta,\gamma\rangle_{ABC}:=[\alpha y-x\gamma]\in H_{A}^{k_1 + k_2 + k_3-2}(X)/([\alpha]H_A(X) + H_A(X)[\gamma]),
    \]
where $\alpha\beta=i\del\delbar x$ and $\beta\gamma=i\del\delbar y$ for some $x\in A^{k_1+k_2-2}$ and $y\in A^{k_2+k_3-2}$.
\end{definition}

\begin{rem}\label{rem: definings systems}\,
\begin{enumerate}
    \item A choice of representatives and primitives $\alpha,\beta,\gamma,x,y$ as in \Cref{def: ABCMP} will be referred to as a \emph{defining system.}
    \item There is in general no natural ring structure on $H_A(X)$. However, the wedge product induces a module structure $H_{BC}\otimes H_A\rightarrow H_A$. In particular, $[\alpha]H_A(X) + H_A(X)[\gamma]$ is a well-defined submodule of $H_A(X)$, called the \textit{indeterminacy}.
    \item Different choices of defining systems in the construction of the Massey product can lead to different elements of $H_{A}(X)$, but the result in the quotient by the indeterminacy is independent of the choice of defining system. In fact, as a set, the coset $\langle\alpha,\beta,\gamma\rangle_{ABC}$ is exactly the set of all classes obtained from all possible defining systems.
    \item   In particular, $\langle \alpha,\beta,\gamma\rangle_{ABC}$ depends only on the cohomology classes $[\alpha],[\beta],[\gamma]$. To keep notation slim, one often still writes $\langle \alpha,\beta,\gamma\rangle_{ABC}$ instead of $\langle [\alpha],[\beta],[\gamma]\rangle_{ABC}$ (as we did above).  
    \item A Massey product of three cohomology classes is said to be trivial, if it is the zero class in the quotient of $H_{A}(X)$ modulo indeterminacy. Equivalently, it is trivial if there exists a defining system that produces the zero class in $H_{A}(X)$.
    \end{enumerate}  
\end{rem}

Our conventions differ slightly from the references given in order to have less signs and be defined over the reals in the sense that if all three inputs are conjugation invariant classes, then so is the output. The $ABC$-Massey products are invariants of the pluripotential quasi-isomorphism type of $A(X)$, \cite[Proposition 4.4]{MilSt_bigrform}.

There are natural maps 
\begin{equation}\label{eqn: comparision a-dr-bc}
    H_{BC}(X)\to H_{dR}(X)\to H_A(X)
\end{equation}
induced by the identity. One says that $X$ satisfies the $\del\delbar$-Lemma if the natural maps \eqref{eqn: comparision a-dr-bc} are isomorphisms. As is well-known \cite{DGMS}, compact K\"ahler (and hence projective) manifolds satisfy the $\del\delbar$-Lemma. For this reason, we will mostly not distinguish these three cohomologies and instead simply write $H(X)$.

We will frequently use the following technique to show an ABC-Massey product is nontrivial, which for ordinary Massey products goes back (at least) to \cite{Taylor_IMP}: Given classes $[\alpha],[\beta],[\gamma]$ such that $\langle \alpha,\beta,\gamma\rangle_{ABC}$ is defined and $[\delta]\in H_{BC}(X)$ such that $[\alpha][\delta]=[\gamma][\delta]=0\in H_{BC}(X)$, then multiplication with $[\delta]$ annihilates the indeterminacy and so the product $\langle\alpha,\beta,\gamma\rangle_{ABC}[\delta]\in H_A(X)$ is a well-defined class. In particular, if $X$ is compact and this product is in top degree, we can integrate it over the manifold. If the result is nonzero, the Massey product was not trivial.

While it will never be needed in the following, we now recall the definition of classical Massey products for context.
\begin{definition}[\cite{Massey_mexico},\cite{Massey_holn}]
    Let $X$ be a manifold and $[\alpha]\in H_{dR}^{k_1}(X)$, $[\beta]\in H_{dR}^{k_2}(X)$, $[\gamma]\in H_{dR}^{k_3}(X)$ s.t. $[\alpha][\beta]=[\beta][\gamma]=0$. The triple Massey product of these classes is the coset
    \[
    \langle \alpha,\beta,\gamma\rangle := \langle (-1)^{|\alpha|}\alpha y - x\gamma\rangle\in H_{dR}^{k_1+k_2+k_3-1}(X)/([\alpha]H_{dR}(X)+H_{dR}(X)[\gamma]),
    \]
    where $\alpha\beta=dx$ and $\beta\gamma=dy$ for some $x\in A^{k_1+k_2-1}(X)$ and $y\in A^{k_2+k_3-1}(X)$.
\end{definition}

There is a natural map $H_A\to H_{BC}$ induced by $d^c=\frac{i}{2}(\delbar-\del)$. With this and the natural maps \eqref{eqn: comparision a-dr-bc} understood, the $ABC$-Massey products are secondary variants of the ordinary Massey products in the sense that one has the following relation in $H_{dR}(X)$:
\begin{equation*}
d^c\langle\alpha,\beta,\gamma\rangle_{ABC}=\langle \alpha,\beta,\gamma\rangle,
\end{equation*}
up to a universal sign depending on one's convention for the ordinary triple Massey products. 
There are versions of (ABC and classical) Massey products taking in $n$ classes for any $n\geq 3$, but we will only consider triple products. By \cite{DGMS}, on a compact K\"ahler manifolds, all classical Massey products are trivial.

\section{Riemann surfaces}

Let $\Sigma$ be any compact Riemann surface. Recall (see e.g. \cite{Donaldson}) there is a positive definite sesquilinear form
\begin{align*}
\langle ~,~\rangle_D:A^{1,0}\times A^{1,0}&\longrightarrow \C\\
(\alpha,\beta)&\longmapsto \langle\alpha,\beta\rangle_D:=i\cdot\int_\Sigma\alpha\wedge \overline{\beta}\,.
\end{align*}

\begin{proof}[Proof of \Cref{Thm: RS}]
    Let $\Sigma$ be a Riemann surface of genus at least two and let us fix some point $x\in \Sigma$. Let $\omega_1,\omega_2\in A^{1,0}(\Sigma)$ be two nonzero holomorphic one-forms such that $\langle\omega_1,\omega_2\rangle_D=i\int_\Sigma\omega_1\wedge\bar{\omega}_2=0$. Then, there exists a unique function $F$ s.t. $F(x)=0$ and $i\del\delbar F= \omega_1\wedge\bar{\omega}_2$. Note that restricted to some simply connected dense open subset $U\subset \Sigma$, we may write $\omega_i|_U=f_idz$ for some holomorphic function $f_i$ and so $\omega_1\wedge\bar{\omega}_2|_U=f_1\bar{f}_2dzd\bar{z}$ which can only have a discrete set of zeroes by the identity principle. 
    
   We first observe that $f_1\bar{f_2}$ is not a real function. In fact, otherwise, since the complement of the zeroes of $f_1\bar{f}_2$ is connected, we have $f_1\bar{f}_2\geq 0$ or $f_1\bar{f}_2\leq 0$. Then, denoting $U'\subseteq U$ the complement of the zeroes of $f_1\bar{f}_2$, we have 
   \[
    0=\langle \omega_1,\omega_2\rangle_D=\int_{U'} f_1\bar{f}_2 idzd\bar{z}\neq 0
   \]
   yielding a contradiction.

   We now set  $\alpha_i=\omega_i+\bar{\omega}_i$. Then, restricted to $U$, the form $\alpha_1\wedge\alpha_2$ has the expression $2i\mbox{Im}(f_1\bar{f_2})dzd\bar{z}$, which is nonzero by assumption and therefore there exists a nonconstant real function $u$ with $i\del\delbar u= \alpha_1\wedge \alpha_2$ (in fact, up to a constant, we may take $\mbox{Re}(F)$). Then, we have
    \[
    0\neq \langle\alpha_1,\alpha_1,\alpha_2\rangle_{ABC}=[u\alpha_1]\in H^1(\Sigma)/([\alpha_1],[\alpha_2])
    \]
    In fact, since
    $[\alpha_1\wedge \alpha_2]=0=[\alpha_2\wedge\alpha_2]$, this follows from
    \[
    \int_\Sigma u\alpha_1\wedge\alpha_2= \int_\Sigma u\cdot i\del\delbar u=-i\cdot \int_\Sigma\del u\wedge\delbar u=-\langle \del u,\del u\rangle_D\neq 0
    \]
    using Stokes' theorem and that $u$ is real.
\end{proof}

\begin{cor}\label{cor: fibrations}
    Any compact K\"ahler manifold which admits a nonconstant holomorphic map to a closed Riemann surface of genus $g\geq 2$ admits a nontrivial $ABC$-Massey product.
\end{cor}
\begin{proof}
    Given such $\pi:Y\to X=\Sigma_{g\geq 2}$, the induced map in cohomology is injective \cite{Wells}. Let $\mathfrak m =\langle \alpha,\beta,\gamma\rangle_{ABC}\in H^{1}(X)/([\alpha],[\gamma])$ be a nonzero Massey product. Then, we have $\pi^*\mathfrak m=\langle\pi^*\alpha,\pi^*\beta,\pi^*\gamma\rangle_{ABC}\neq 0\in H^1(Y)/(\pi^*[\alpha],\pi^*[\gamma])$, where we use that in degree $1$ the ideal generated by $[\pi^*\alpha], [\pi^*\gamma]$ equals the image under $\pi^*$ of the ideal generated by $[\alpha]$ and $[\gamma]$, i.e. the relevant indeterminacy in $Y$ is not bigger than in $X$.
\end{proof}
\begin{rem}\label{rem: rigid}
    Note that this means there exist rigid manifolds which are not pluripotentially formal, as there are rigid manifolds fibering over curves of genus $\geq 2$.
\end{rem}
\begin{rem}
    The property of admitting a fibration over a curve of genus $g\geq 2$ depends only on the de Rham cohomology ring \cite{Catanese1991}. So \Cref{cor: fibrations} can be rephrased as follows: No strongly formal K\"ahler structure can exist on a manifold $X$ which admits an isotropic subspace $U\subset H^1(X)$ with $\dim(U)\geq2$.
\end{rem}

\begin{rem}
    The argument in the proof of \Cref{cor: fibrations} works for any holomorphic map $\pi:Y\to \Sigma_{\geq 2}$ without the K\"ahler assumption on $Y$ as long as $\pi^*$ is injective in $H^1$. What is perhaps interesting is that in that case the pullback of the class $\langle\alpha_1,\alpha_1,\alpha_2\rangle_{ABC}[\alpha_2]$ (which we used to certify non-triviality of the product) can indeed be trivial: This happens for example if one takes $Y$ to be the quotient of the form $(L\backslash\{0\})/\langle\lambda^\Z\rangle$, where $L\to \Sigma$ is a holomorphic line bundle with nonzero first Chern class over $\Sigma$ and $\lambda$ is a nonzero complex number of modulus $\|\lambda\|\neq 1$ acting in the fibres: In fact, in this case the fundamental class of $\Sigma$ is a nonzero multiple of the first Chern class of $L$ and hence pulls back trivially.
\end{rem}

\section{A simply-connected example}

Let $P,Q,R\in\CP^3$ be in general position and let $C$ be the line through $P,Q$. Choose distinct points $S_i\in C\backslash\{P,Q\}$ for $i=1,2$ and let $L_i$ be the line through $R$ and $S_i$. Now let $X$ be the space obtained from $\CP^3$ by first blowing up $P,Q,R$, then blowing up the strict transform of $C$, and finally blowing up the strict transforms of the $L_i$. We denote by $E_P,E_Q,\ldots$ the (strict transforms of the) respective exceptional divisors in $X$. Let furthermore $A_i\subset X$ be the strict transform of a hyperplane containing $L_i$ and intersecting $C$ transversally. 


\begin{minipage}{6in}
  \centering
  $\vcenter{\hbox{\includegraphics[scale=0.3]{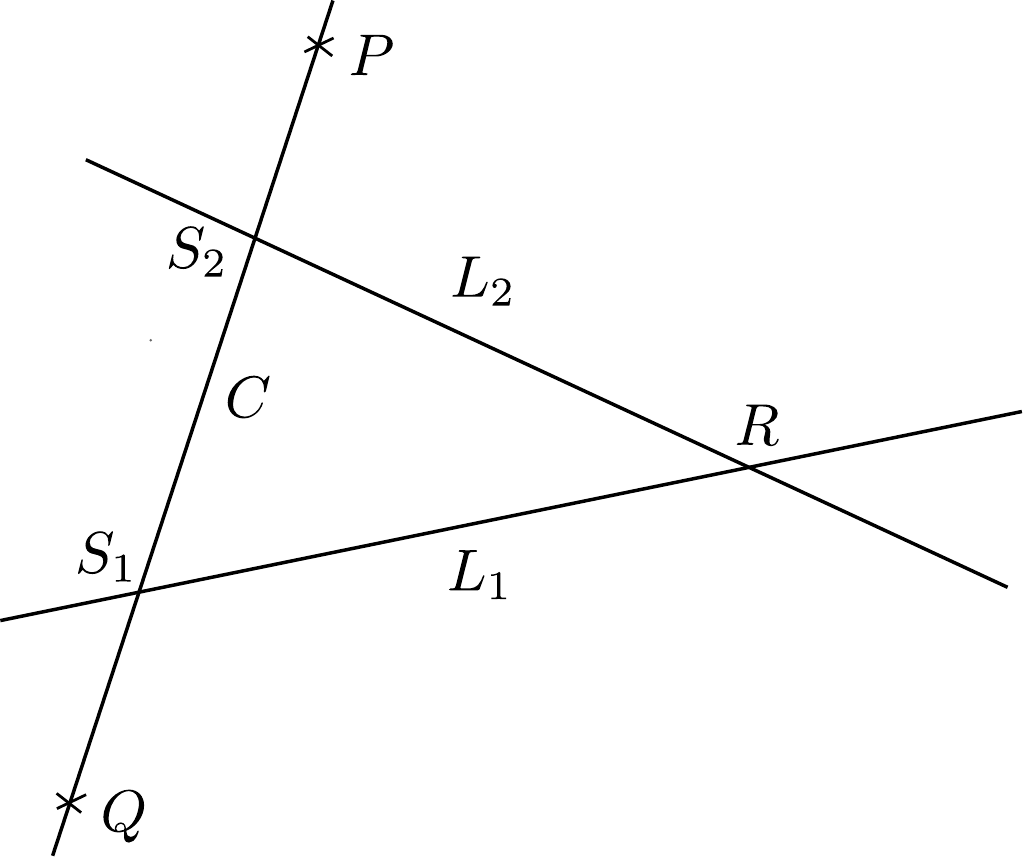}}}$
  \hspace*{.3in} $\longsquiggly$ \hspace*{.3in}
  $\vcenter{\hbox{\includegraphics[scale=0.3]{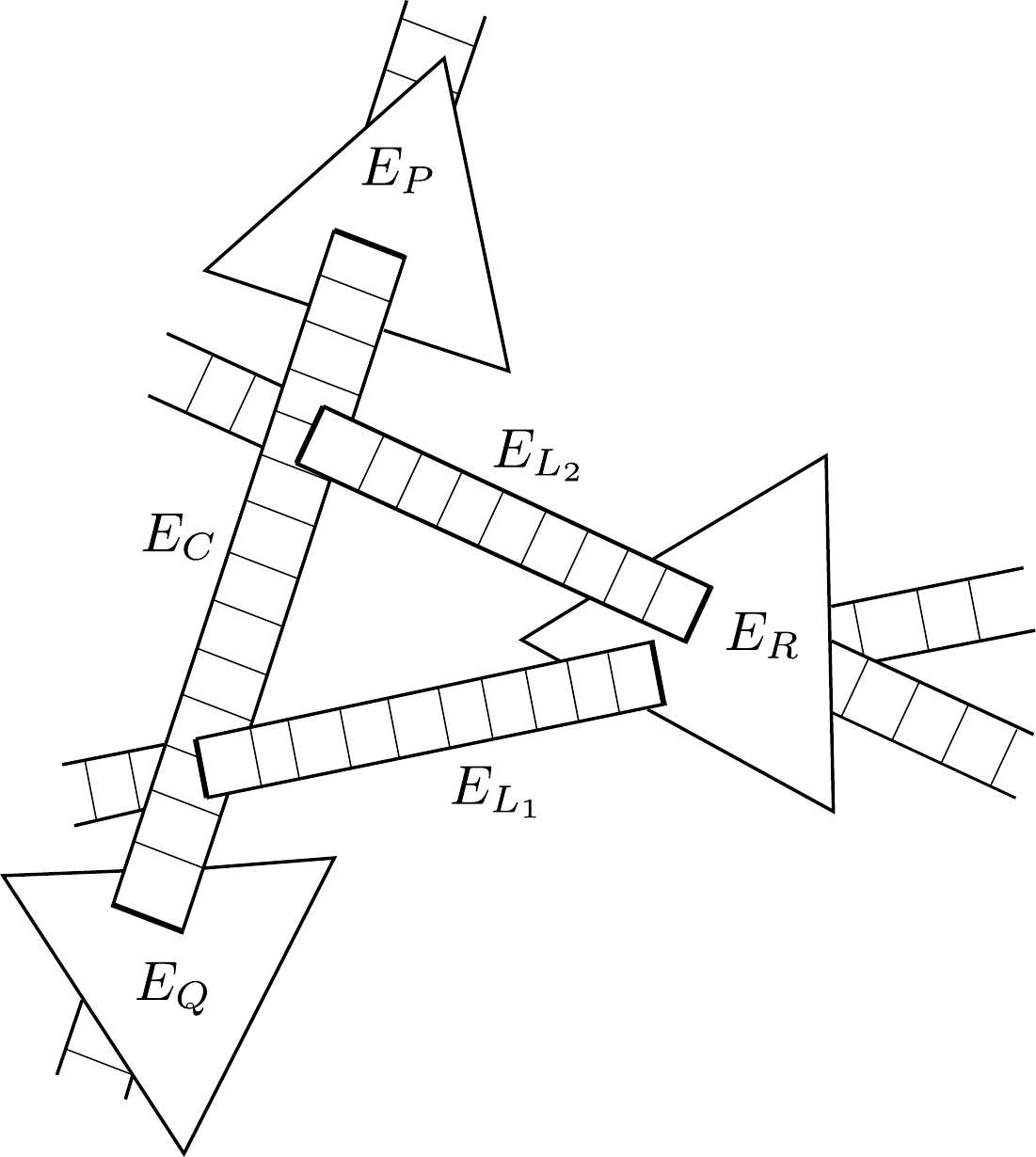}}}$
  \end{minipage}

Recall the definition of the cross-ratio of four distinct points $P_i\in \CP^1=\C\cup\{\infty\}$:
\begin{equation}
    \chi(P_1,P_2,P_3,P_4):=\frac{(P_3-P_1)(P_4-P_2)}{(P_3-P_2)(P_4-P_1)},
\end{equation}
where one applies the naive calculation rules in case one of the points is $\infty$, e.g. $0=\frac{1}{\infty}$. This quantity is invariant under projective transformations and hence well defined for any $4$-tuple of distinct points on a line. The main objective of this section is to prove

\begin{thm}\label{thm: mainthm simplycon}
    The map \begin{align*}
        H^{2,2}(X)&\longrightarrow \C\\
        x&\longmapsto \int_{A_1-A_2}x
    \end{align*}vanishes on the indeterminacy of the Massey product
\[\mathfrak m:=\langle [E_{L_1}+E_{L_2}+E_R],[E_P-E_Q],[E_C-E_R]\rangle_{ABC}\]
and 
\[\int_{A_1-A_2}\mathfrak{m}=
\frac{\log|\chi(P,Q,S_1,S_2)|}{\pi}.
\]
\end{thm}

Let us first check that the Massey product $\mathfrak m$ is indeed defined.  We use that for cycles that meet transversely, the cup product of their fundamental classes is dual to the fundamental class of their intersection. The divisors $E_{L_1}+E_{L_2}+E_R$ and $E_P-E_Q$ are disjoint so the product on the left hand side of the Massey product is indeed $0$. On the right hand side, we have $[E_P-E_Q]\cdot [E_C-E_R]=[E_C\cap E_P-E_C\cap E_Q]$, which is the trivial class since all fibres of $E_C\to C$ are rationally, and hence homologically, equivalent. 
To exhibit an explicit primitive for $\delta_{E_C\cap E_P-E_C\cap E_Q}$, let us identify $C\cong \CP^1=\C\cup\{\infty\}$ in such a way that $P=0$ and $Q=\infty$. We may then view the projection $f:E_C\to C$ as a meromorphic function. Thus, by \Cref{prop:LLP} applied to $f$ and the trivial line bundle, we have
\[\delta_{E_C\cap E_P}-\delta_{E_C\cap E_Q}= \frac i \pi \del\delbar \log |f|\]
as currents on $E_C$ and the push-forward along the inclusion $j:E_C\rightarrow X$ yields the desired $i\del\delbar$-primitive $p:=j_*(\frac 1 \pi \log|f|)\in\mathcal{D}^{1,1}(X)$ certifying triviality of the product of the right hand cohomology classes in $\mathfrak m$.

\begin{rem}
  Let us now do a heuristic calculation of the Massey product  using the "defining system" \[(\delta_{E_{L_1}+E_{L_2}+E_R}, \delta_{E_P-E_Q}, \delta_{E_C-E_R}, 0, p)\] which consists of non-smooth currents and using the partial multiplicative structure from \Cref{rem: product currents}. We arrive at the following representative for $\mathfrak{m}$
    \[
    \frac{1}{\pi}\delta_{E_{L_1}+E_{L_2}+E_R}\cdot p =     \frac{1}{\pi}j_*(\log|f|\cdot \delta_{F_1} + \log |f|\cdot \delta_{F_2}) = \frac 1 \pi (\log|S_1|\cdot \delta_{F_1}+\log|S_2|\cdot \delta_{F_2})~,
    \]
    where $F_i=E_C\cap E_{L_i}$ and $S_1,S_2$ are considered points in $\CP^1$ with respect to the above identification $C\cong \CP^1$. The first equality is the projection formula and the second follows because $f|_{F_i}$ is constant.
    
   However, while the inclusion of forms to currents is a quasi-isomorphism, since the multiplication on currents is only partially defined, it is not a priori clear that the above indeed represents a class in $\mathfrak{m}$. To avoid these technical issues, in the next lemma, we compute $\mathfrak{m}$ rigorously using a defining system consisting of smooth forms.
\end{rem}

\begin{lem}\label{lem: actual work}
Let $F_i=E_C\cap E_{L_i}\subset X$.
With respect to the above coordinates on $C$ one has $\frac 1 \pi (\log|S_1|[F_1]+\log|S_2|[F_2])\in\mathfrak m$.
\end{lem}

\begin{proof}
    We set $D_1=E_{L_1}+E_{L_2}+E_R$, $D_2=E_P-E_Q$, and $D_3=E_C-E_R$. In order to compute $\mathfrak m=\langle [D_1],[D_2],[D_3]\rangle_{ABC}$ we need to find a defining system for $\mathfrak m$ in $A(X)$.
    
    As a first step we choose suitable differential form representatives for the cohomology classes $[D_i]$. Let $\Oh(D_i)$ be the line bundle with section $s_i$ corresponding to the divisor $D_i$. Then for any Hermitian metric $h_i$ on $\Oh(D_i)$ we obtain the Chern form $\gamma_i:=c_1^{h_i}(D_i)$. 

     Recall from \Cref{prop:LLP} that
    \[
        \frac i \pi \del\delbar \log |s_i|_{h_i}=\delta_{D_i}-\gamma_i
    \]
    so we may take $\gamma_i$ as a form representative for $[D_i]$.
    Since $\Oh(D_i)$ is trivial outside $D_i$, one may always pick $h_i$ such that $|s_i|\equiv 1$ except on some small neighborhood of $D_i$. Hence we may choose $\log |s_i|_{h_i}$ and in particular $\gamma_i$ to be supported near $D_i$.

    To complete the defining system and hence compute the value of the Massey product we now need to find explicit $i\del\delbar$-primitives for the exact forms $\gamma_1\gamma_2$ and $\gamma_2\gamma_3$. Due to the supports we have $\gamma_1\gamma_2=0$ and may choose $0$ as a $i\del\delbar$-primitive for $\gamma_1\gamma_2$. Let $p=j_*(\frac 1 \pi \log|f|)$ be the previously constructed $i\del\delbar$-primitive for the current $\delta_{D_2\cap D_3}$. We use $p$ to find a form primitive for $\gamma_2\gamma_3$. First, set
    \[
    q:=p-\frac{1}{\pi}\log|s_2|_{h_2} \delta_{D_3}-\frac{1}{\pi}\log|s_3|_{h_3} \gamma_{2}\in \mathcal{D}^{1,1}(X),
    \] 
    where we note that the middle term is a well-defined current since $\log|s_2|_{h_2}$ restricts to a locally integrable function on $D_3$.
     We have

     \begin{align*}
        i\del\delbar q&= \delta_{D_2\cap D_3}-\frac{i}{\pi}\del\delbar \left(\log |s_2|_{h_2} \delta_{D_3}\right)-\frac{i}{\pi}\del\delbar\left(\log|s_3|_{h_3} \right)\gamma_{2}\\
        &= \delta_{D_2\cap D_3}-\left(\delta_{D_2\cap D_3}-\gamma_2\delta_{D_3}\right)-\left(\delta_{D_3}-\gamma_3\right) \gamma_{2}\\
        &=\gamma_2\gamma_3~,
     \end{align*}
    where for the first equality we use the Leibniz rule on the right hand summand and for the second equality we apply \Cref{prop:LLP} (note that on the middle summand it is applied to respective restrictions of $\log|s_2|_{h_2}$ with the manifolds in $D_3$ serving as the ambient manifold).

   Now, if we choose a smooth $i\del\delbar$-primitive $r\in A^{1,1}(X)$ for $\gamma_2\gamma_3$, the difference $r-q$ defines an Aeppli class. By adding a smooth representative for the negative of this class, we may assume that $r-q$ defines a trivial class in $H_A(X)$. We take such $r$ as a primitive for $\gamma_2\gamma_3$ in our defining system for $\mathfrak{m}$. Thus, with the defining system $\gamma_1,\gamma_2,\gamma_3, 0, r$ the ABC-Massey product $\mathfrak m$  is then represented by $[\gamma_1 r]\in H^{2,2}_A(X)$.

    It remains to see that this cohomology class agrees with the one in the statement of the lemma. If one allows current representatives of cohomology classes, one in fact has
    \[
    [\gamma_1r]=[\gamma_1 q]=[\gamma_1 p]
    \]
    where the first equality is due to $[r-q]=0$ and for the second equality we use that $p-q$ is supported near $D_2\cap D_3$ and hence has support disjoint from the support of $\gamma_1$. Now by the definition of $p$ we have $\gamma_1 p= j_*(\frac 1 \pi \log|f|\cdot \gamma_1|_{E_C})$.
    Note that $\log|f|\cdot \gamma_1|_{E_C}$ is supported inside a small open neighbourhood $U\subset E_C$ of $D_1\cap E_C$ for which $[\gamma_1|_{U}]$ and $[D_1\cap E_C]$ agree in the Bott-Chern cohomology $H_{BC,cpt}(U)$ with compact support by \Cref{prop:LLP}. Furthermore the restriction of $\log|f|$ to $U$ defines an element in $H_A(U)$ and, using that the pairing
    \[H_A(U)\otimes H_{BC,cpt}(U)\rightarrow H_{A,cpt}(U)\]
    is well-defined, we find that $[\log |f|\cdot \gamma_1|_{E_C}]$ and $[\log|f| \cdot \delta_{D_1\cap E_C}]$ agree in $H_{A,cpt}(U)$. But then so do their images in $H_A(E_C)$ under the pushforward of the inclusion $U\rightarrow E_C$. This shows that $\mathfrak{m}$ is represented by $\frac 1 \pi \log|f|\cdot \delta_{D_1\cap E_C}=\frac 1 \pi \log|f|\cdot \delta_{F_1+F_2}=\frac{1}{\pi}(\log|f|\cdot \delta_{F_1}+\log|f|\cdot \delta_{F_2})$. By its definition $f$ is constant on $F_i$ and takes the value $S_i$ with respect to the initial identification $C=\C\cup\{\infty\}$.
\end{proof}

\begin{rem}
    The calculations in \Cref{lem: actual work} happen away from $E_R$ and equally show that
    \[\frac 1 \pi (\log|S_1|\cdot [F_1]+\log|S_2|\cdot [F_2])\in \langle [E_{L_1}+E_{L_2}],[ E_P-E_Q], [E_C]\rangle_{ABC}.\]
    However, while we will show that $\mathfrak m$ is nontrivial, the simplified Massey product above is in fact trivial: the $[F_i]=[E_C]\cdot[E_{L_i}]$ lie in the indeterminacy and hence the Massey product represents the trivial coset.
\end{rem}

To deal with the indeterminacy of $\mathfrak m$, we will now compute the multiplicative structure of $H(X)$. Recall the formulae for the blow-up of a submanifold $Z\subseteq Y$ of codimension $r$:
    \[
    H(Bl_Z Y)\cong H(Y)\oplus \bigoplus_{i=1}^{r-1} H(Z)E_Z^i
    \]
 The algebra structure on the right is given as follows (see e.g. \cite[Thm. 7.31]{Voisin}): $H(Y)$ is a subring, $H(Z)$ is an $H(Y)$-module via restriction, and $E_Z^{r}=(-1)^{r-1} [Z] + \sum_{i=1}^{r-1} (-1)^{r-i+1}c_i(\mathcal{N}_{Z})E^{r-i}_Z$, where $\mathcal{N}_{Z}$ is the normal bundle. The isomorphism from right to left is given by pullback on $H(Y)$, sends the symbol $E_Z$ to $[E_Z]$, and is a map of $H(Y)$-modules on $H(Z)$.

 It follows that $H^{1,1}(X)$ is generated by the pullbacks of the classes of a generic plane $G$ in $\CP^3$ as well as the exceptional divisors $E_P$, $E_Q$, $E_R$, $E_C$, $E_{L_1}$, and $E_{L_2}$. Note that the pullback of the Thom class of a divisor along a blow down is the Thom class of its strict transform, provided the centre of the blow-up is not contained in the divisor. For ease of notation, in the following we omit square brackets around these divisors when indicating the corresponding cohomology class in $H(X)$. One obtains that a basis for $H^{2,2}(X)$ is given by $G^2$, $E_P^2$, $E_Q^2$, $E_R^2$, $GE_C$, $GE_{L_1}$, and $GE_{L_2}$, where we have used that $G$ restricts to a generator of $H^{1,1}(C)$, $H^{1,1}(L_1)$, and $H^{1,1}(L_2)$.

 \begin{rem}\label{rem: geometric rep's}
    Geometric representatives for this basis for $H^{2,2}(X)$ are given as follows:
    \begin{itemize}
        \item $G^2$ equals the class of a generic line (not meeting any exceptional set).
        \item $E_P^2, E_Q^2$, resp.  $E_R^2$ equal  minus the classes of a generic line in $E_P, E_Q$, resp. $E_R$.
        \item $GE_C, GE_{L_1}, GE_{L_2}$ equal the classes of a general fibre of $E_C\to C$, resp. $E_{L_i}\to L_i$.
    \end{itemize}
 \end{rem}
 The following table displays the map $H^{1,1}(X)\otimes H^{1,1}(X)\rightarrow H^{2,2}(X)$ in the chosen bases. 
\vspace{0.3cm}\begin{figure}[h]
    \centering
\begin{tabular}{c||c|c|c|c|c|c|c}
    &$G$&$E_P$&$E_Q$&$E_R$&$E_C$&$E_{L_1}$&$E_{L_2}$\\
    \hline\hline
    $G$&$G^2$&$0$&$0$&$0$&$GE_C$&$GE_{L_1}$&$GE_{L_2}$\\
    \hline
    $E_P$&0&$E_P^2$&$0$&$0$&$GE_C$&$0$&$0$\\
    \hline
    $E_Q$&0&0&$E_Q^2$&$0$&$GE_C$&$0$&$0$\\
    \hline
    $E_R$&0&0&0&$E_R^2$&$0$&$GE_{L_1}$&$GE_{L_2}$\\
    \hline
    $E_C$&$GE_C$&$GE_C$&$GE_C$&0&\makecell{$-2GE_C$\\$-G^2-E_P^2-E_Q^2$}&$GE_{L_1}$&$GE_{L_2}$\\
    \hline
    $E_{L_1}$&$GE_{L_1}$&$0$&$0$&$GE_{L_1}$&$GE_{L_1}$&\makecell{$-GE_{L_1}$\\$-G^2-E_R^2+G E_C$}&$0$\\
    \hline
    $E_{L_2}$&$GE_{L_2}$&0&0&$GE_{L_2}$&$GE_{L_2}$&0&\makecell{$-GE_{L_2}$\\$-G^2-E_R^2+G E_C$}
\end{tabular}

    \caption{Multiplication table for $H(X)$.}
    \label{fig:m-table}
\end{figure}

The off-diagonal entries can be computed via the transverse intersections of the divisors and the fact that all fibers in the blowup of a line are rationally equivalent.  Regarding the diagonal entries, only $E_C^2$, $E_{L_1}^2$, and $E_{L_2}^2$ need justification. We provide the computation for $E_C^2$. The computations for $E_{L_i}$ are similar and left to the reader. In the above formula for the blowup cohomology we get $E_C^2=-C+ c_1(\mathcal{N}_C) E_C$, where the right hand side is to be understood in the blow-up of $\CP^3$ at the points $P,Q,R$, say $X_1$, where we are about to blow up $C$. We determine the class $C$ in $H^{2,2}(X_1)$ by pairing it geometrically with the basis $G,E_P,E_Q,E_R$ of $H^{1,1}(X_1)$. We find that $C=G^2+E_P^2+E_Q^2$. We further compute \[\int_Cc_1(\mathcal{N}_C)=\int_C(c_1(X_1)|_C-c_1(C))=\int_C(4G-2E_P-2E_Q-2E_R)|_C-2=-2.\]

\begin{lem}\label{lem:indeterminacy}
    Let $I\subset H^{*,*} (X)$ be the submodule generated by $E_C-E_R$ and $E_R+E_{L_1}+E_{L_2}$. Then for any representative of an element of $I^{2,2}$ the integrals over $A_1$ and $A_2$ agree.
\end{lem}

\begin{proof}
 From the multiplication table above one computes that \[I^{2,2}=\langle GE_C,E_R^2,E_P^2+E_Q^2,G(E_{L_1}+E_{L_2}), G^2\rangle.\] 
 All intersection numbers between the $A_i$ and the generating cycles of $I^{2,2}$ can be computed directly via transverse intersections with geometric representatives listed in \Cref{rem: geometric rep's} and agree for $i=1,2$.
    \end{proof}

    \begin{rem}
  The following slightly stronger statement holds true: The class $[A_1-A_2]$ pairs trivially with the generators of the indeterminacy submodule for $\mathfrak{m}$.
  Indeed, one can compute that the class $[A_i]$ agrees with $G-E_{L_i}-E_{R}$. The pairing with the generators of the ideal can then be computed from the multiplication table.
\end{rem}

\begin{proof}[Proof of Theorem \ref{thm: mainthm simplycon}.] By Lemma \ref{lem:indeterminacy}, multiplying with $[A_1-A_2]$ sends the indeterminacy ideal $I$ to $0$. From the computation in Lemma \ref{lem: actual work} and since $F_i$ and $A_i$ meet transversally in a point, we deduce that \[
\pi \cdot \int_{A_1-A_2}\mathfrak m=\log|S_1|-\log|S_2|=\log\left|\frac{S_1}{S_2}\right|.\] As $P=0$ and $Q=\infty$ in the chosen coordinate on $C$, this is indeed $\log|\chi(P,Q,S_1,S_2)|$.
 \end{proof}

As an immediate consequence of Theorem \ref{thm: mainthm simplycon} we obtain

\begin{cor}\label{Cor: blowups nonform}
    For any complex manifold $Y$ of dimension at least $4$ there is a finite sequence of blow-ups of points and lines such that the resulting manifold $\widetilde Y$ admits a non trivial triple $ABC$-Massey product and is thus in particular not strongly formal.
\end{cor}

\begin{proof}
    Let $Y$ be a complex manifold of dimension $n\geq 4$. After possibly blowing up a point in $Y$, we can assume the existence of an embedding $\CP^3\rightarrow Y$. Now choose a configuration of points $P,Q,R$ and lines $C,L_i$ inside $\CP^3$ and blow up $Y$ along these submanifolds as in the construction of $X$. The Massey product 
    \[\langle [E_{L_1}+E_{L_2}+E_R],[E_P-E_Q],[E_C-E_R]\rangle_{ABC}\] is defined in the resulting space $\widetilde Y$. Furthermore there is an embedding $X\rightarrow\widetilde Y$ given by the strict transform of $\CP^3$. Cohomologically it maps the class of the (strict transform of the) exceptional divisors to their respective intersections with $X$, which are the (strict transforms of the) exceptional divisors used to define $\mathfrak m$ in Theorem \ref{thm: mainthm simplycon}. Hence there is an $ABC$-Massey product $\widetilde{\mathfrak m}\subset H_A(\widetilde Y)$ which pulls back to a subset of $\mathfrak m\subset H_A(X)$ and is thus nontrivial for suitable choices in the construction.
\end{proof}

\begin{cor}\label{biholotypes}
    For uncountably many choices of $P,Q,R,S_1,S_2$, the resulting manifolds $X$ as in \Cref{thm: mainthm simplycon} are pairwise non-biholomorphic, and distinguished by the $ABC$-Massey product. 
\end{cor}
\begin{remark}\label{rem: no odd}
Note that the examples in \Cref{biholotypes} have no odd cohomology. In particular, they cannot be distinguished using intermediate Jacobians.
\end{remark}
\begin{proof}
    Consider a biholomorphism $\varphi:X'\to X$ between two manifolds corresponding to different choices, say $X=X(P,Q,R,S_1,S_2)$ and $X'=X(P',Q',R',S_1',S_2')$, and the corresponding pullback $H(X)\to H(X')$. In $X'$, there is at most a countable collection of Massey products with integral $(1,1)$-classes as entries and correspondingly only a countable collection of values when integrating these Massey products over integral cycles annihilating the indeterminacy. Thus, for any such $\varphi$, $\int_{A_1-A_2}\mathfrak m = \int_{\varphi^{-1}(A_1-A_2)}\varphi^*\mathfrak m$ must be among these values.
\end{proof}
\begin{rem}
    One may be more precise and show that the automorphism group of $H(X)$ is generated by the involutions exchanging $E_P$ and $E_Q$, resp. $E_{L_1}$ and $E_{L_2}$. Thus, whenever the values $|\log|\chi(P,Q,S_1,S_2)||$ and $|\log|\chi(P',Q',S_1',S_2')||$ differ, the manifolds $X$ and $X'$ are not biholomorphic.
\end{rem}

    \bibliographystyle{amsplain}

\bibliography{bib.bib}

\end{document}